\documentclass{amsart}

\usepackage{amsmath,amssymb,amsfonts,enumerate,amsthm}
\usepackage[all]{xy}

\newtheorem{lem}{Lemma} 
\newtheorem{cor}[lem]{Corollary}

\newtheorem{thm}[lem]{Theorem}

\newtheorem{Defn}[lem]{Definition}
\newtheorem{Ex}[lem]{Example}
\newtheorem{Question}[lem]{Question}
\newtheorem{Property}[lem]{Property}
\newtheorem{Properties}[lem]{Properties}
\newtheorem{Discussion}[lem]{Discussion}
\newtheorem{Construction}[lem]{Construction}
\newtheorem{Subprops}{}[lem]
\newtheorem{Para}[lem]{}

\newcommand{\Ann}{\operatorname{Ann}}
\newcommand{\Hom}{\operatorname{Hom}}
\newcommand{\Ext}{\operatorname{Ext}}
\newcommand{\Tor}{\operatorname{Tor}}

\newcommand{\Ass}{\operatorname{Ass}}

\newcommand{\depth}{\operatorname{depth}}

\newcommand{\pd}{\operatorname{pd}}
\newcommand{\id}{\operatorname{id}}

\renewcommand{\H}{\operatorname{H}}

\newcommand{\G}{\operatorname{G}}

\newcommand{\fa}{\mathfrak{a}}

\newcommand{\fm}{\mathfrak{m}}
\newcommand{\fp}{\mathfrak{p}}

\newcommand{\gkdim}[1]{\mathrm{G}_{#1}\text{-}\dim}     
\newcommand{\DD}{\mathcal{D}}
\newcommand{\xra}{\xrightarrow}
\newcommand{\zz}{\mathbb{Z}}
\newcommand{\rhom}{\mathbf{R}\mathrm{Hom}}      
\newcommand{\lotimes}{\otimes^{\mathbf{L}}}
\newcommand{\gdim}{\mathrm{G\text{-}dim}}       
\newcommand{\shift}{\Sigma}
\newcommand{\ol}{\overline}
\newcommand{\vf}{\varphi}
\newcommand{\lcidim}{\mathrm{CI}_{\ast}\text{-}\mathrm{dim}}
\newcommand{\gid}{\operatorname{Gid}}
\newcommand{\coker}{\operatorname{Coker}}
\newcommand{\HH}{\operatorname{H}}

\begin{document}

\bibliographystyle{amsplain}

\author{Shokrollah Salarian}
\address{Department of Mathematics, Isfahan University, P.O. Box 81746-73441, 
Isfahan, Iran, and Institute for Theoretical Physics and Mathematics}
%Institute for Studies in Theoretical Physics and
%Mathematics, and Department of Mathematics, Damghan University}
\email{salarian@ipm.ir}

\author{Sean Sather-Wagstaff}
\address{Department of Mathematics, University of Nebraska, 
203 Avery Hall, Lincoln, NE, 68588-0130 USA}
\email{swagstaff@math.unl.edu}

\author{Siamak Yassemi}
\address{Department of Mathematics, University of Tehran, and 
Institute for Studies in Theoretical Physics and Mathematics}
\email{yassemi@ipm.ir}

%\thanks{\today}
\thanks{The research of the first author 
was supported in part by a grant from IPM (no. 84130031). This research 
was conducted while the second author was an NSF Mathematical Sciences 
Postdoctoral Research Fellow.  The research of the third author was 
supported in part by a grant from the IPM (no. 82130212).}



\title[Characterizing local rings]{Characterizing 
local rings via homological dimensions and regular 
sequences}

\begin{abstract}
Let $(R,\fm)$ be a Noetherian local ring 
of depth $d$ and $C$ a semidualizing $R$-complex.
Let $M$ be a finite $R$-module and $t$ an integer 
between 0 and $d$.  If $G_C$-dimension of $M/\fa M$ is
finite for all ideals $\fa$ generated by an $R$-regular
sequence of length at most $d-t$ then either
$G_C$-dimension of $M$ is at most $t$ or $C$ is a
dualizing complex.  Analogous results for other homological dimensions 
are also given.
\end{abstract}

\maketitle

\section*{Introduction}

Throughout this work, $(R,\fm, k)$ is a commutative Noetherian local
ring with identity and all modules are unitary.

The characterization of a local ring by the homological properties of 
its finite modules begins with the classical theorem of Auslander, 
Buchsbaum, and Serre, c.f., \cite[(2.2.7)]{bruns:cmr}.  Here $\pd$ 
and $\id$ denote projective and injective dimensions, 
respectively.

\medskip

\noindent \textbf{Theorem A.}  
\textit{The following conditions are equivalent:
\begin{enumerate}[\quad\rm(i)]
\item  $R$ is regular;
\item  $\pd_Rk<\infty$;
\item  $\pd_RM<\infty$ for all finite $R$-modules $M$;
\item  $\id_Rk<\infty$;
\item  $\id_RM<\infty$ for all finite $R$-modules $M$.
\end{enumerate}
}

\medskip

Other criteria for the regularity of $R$ have been given since the 
appearance of this result.  
For example, Foxby~\cite{foxby:dhfmoccmr} offers the following.

\medskip

\noindent \textbf{Theorem A1.}  
\textit{If $M$ is a finite $R$-module such that
$\pd_R(M/\fa M)$ is finite for all ideals $\fa$ of finite
projective dimension, then $M$ is free or $R$ is regular.}

\medskip

Jothilingam and Managayarcarassy~\cite{jothilingam}
are responsible for the next result in this style.  
They prove that the conclusion of Foxby's result holds 
if the hypothesis is satisfied by
ideals generated by regular sequences, the ideal generated
by the empty sequences being
the zero ideal.  This result provides a central 
motivation for our work in this paper.

\medskip

\noindent \textbf{Theorem A2.}  
\textit{If $M$ is a finite $R$-module such that
$\pd_R(M/\fa M)$ is finite for all ideals $\fa$
generated by $R$-regular sequences, then $M$ is free or $R$ is
regular.}

\medskip

Beginning in the late 1960's, several homological dimensions have 
appeared that can be used to detect ring-theoretic properties of 
$R$.  The Gorenstein dimension 
for finite modules was introduced by 
Auslander~\cite{auslander:adgeteac} and developed by Auslander and 
Bridger~\cite{auslander:smt};  see also the monograph of 
Christensen~\cite{christensen:gd}.  
This was extended to $\mathrm{G}_C$-dimension by Foxby~\cite{foxby:gmarm}
and Golod~\cite{golod:gdagpi}, and was studied extensively by
Christensen~\cite{christensen:scatac} and Gerko~\cite{gerko:ohd}.
Other homological dimensions include 
Avramov's virtual projective dimension~\cite{avramov:mofvpd}, 
the complete 
intersection dimension of
Avramov, Gasharov, and Peeva~\cite{avramov:cid},
Veliche's upper 
G-dimension~\cite{veliche:comwfhd}, 
the Gorenstein injective dimension of Enochs and 
Jenda~\cite{enochs1,enochs2},  
and Gerko's lower 
complete intersection dimension 
and Cohen-Macaulay dimension~\cite{gerko:ohd}.  See 
Section~\ref{sec1} for definitions.  

Each of these homological 
dimensions exhibits a theorem characterizing the appropriate 
ring-theoretic property as in Theorem A.  
The goal of this paper is to give analogues of Theorem A for each one.
This is done in Section~\ref{sec2} after the relevant 
background is given in Section~\ref{sec1}.

\section{Homological dimensions} \label{sec1}

We shall employ a small amount of technology from the derived category 
of $R$-modules $\DD(R)$.  We refer the reader to 
Hartshorne~\cite{hartshorne:rad}, Verdier~\cite{verdier:cd}, or 
Gelfand and Manin~\cite{gelfand:moha} for the appropriate background.  

A complex of $R$-modules $X$ is a sequence of $R$-module homomorphisms
\[ \cdots\xra{\partial_{i+1}}X_i \xra{\partial_{i}} X_{i-1}
\xra{\partial_{i-1}}\cdots \]
such that $\partial_{i-1}\partial_i=0$ for each $i\in\mathbb{Z}$.  
The symbol $\simeq$ denotes an isomorphism in $\DD(R)$.
The $i$th homology of $X$ is denoted $\H_i(X)$.  The complex $X$ is 
\emph{homologically bounded} if $\H_i(X)=0$ for almost all
integers $i$;  it is \emph{homologically finite} if it is 
homologically bounded and each $\H_i(X)$ is a finite $R$-module.  
The \emph{infimum of $X$} is 
\[ \inf(X)=\inf\{i\in\zz\mid\H_i(X)\neq 0\}; \]
this is finite when $X\not\simeq 0$ is homologically bounded and is 
$\infty$ for $X\simeq 0$.  

For complexes $X,Y$ the derived tensor product of $X$ and $Y$ is 
denoted $X\lotimes_R Y$ while $\rhom_R(X,Y)$ denotes the derived 
homomorphisms from $X$ to $Y$.  

A homologically finite complex of $R$-modules $C$ is 
\emph{semidualizing} if the natural homothety morphism 
$R\to\rhom_R(C,C)$ is an isomorphism in $\DD(R)$.  When $C$ is a 
finite $R$-module, it is semidualizing if and only if the natural 
homothety morphism $R\to\Hom_R(C,C)$ is an isomorphism and 
$\Ext_R^i(C,C)=0$ for each integer $i\neq 0$.  
The $R$-modules $R$ is always semidualizing.   A semidualizing 
complex is \emph{dualizing} if and only if it has finite injective 
dimension.  

Semidualizing modules 
are studied extensively in~\cite{foxby:gmarm,gerko:ohd,golod:gdagpi}
where they are called ``suitable'' modules.  
Semidualizing complexes are investigated 
in~\cite{christensen:scatac}, whence comes our treatment of 
$\mathrm{G}_C$-dimension, as well as 
in~\cite{frankild:tsosciams,gerko:sdc}.

Let $C$ be a semidualizing $R$-complex.
A homologically finite complex of $R$-modules $X$ is 
\emph{$C$-reflexive} if $\rhom_R(X,C)$ is homologically bounded and 
the natural biduality morphism $X\to\rhom_R(\rhom_R(X,C),C)$ is an 
isomorphism in $\DD(R)$.  The \emph{$\mathrm{G}_C$-dimension of $X$} 
is 
\[ \gkdim{C}_R(X)=\begin{cases}
\inf(C)-\inf(\rhom_R(X,C)) & \text{when $X$ is $C$-reflexive} \\
\infty & \text{otherwise.} \end{cases}
\]
Note that this definition provides $\gkdim{C}_R(0)=-\infty$.
When $C$ and $X$ are both modules, we have the 
following alternate description 
of $\mathrm{G}_C$-dimension of $X$ in terms of resolutions by 
appropriate modules.  

The \emph{$G_C$-class}, denoted $\G_C(R)$, is the collection of 
finite $R$-modules $M$ such that
\begin{enumerate}[\quad(1)]
\item $\Ext^{\ell}_R(M,C)=0$ for all $\ell>0$;
\item $\Ext^{\ell}_R(\Hom_R(M,C),C)=0$ for $\ell>0$; and
\item the biduality morphism 
$M\to\Hom_R(\Hom_R(M,C),C)$ is an isomorphism.
\end{enumerate}
The $\mathrm{G}_C$-dimension of a nonzero finite $R$-module 
$X$ is then the infimum of the set of all nonnegative integers $r$ 
such that there exists an exact sequence 
\[ 0\to M_r\to M_{r-1}\to\cdots\to M_0\to X\to 0 \]
with each $M_i$ in $G_C(R)$.  In particular, $\gkdim{C}_R(M)\geq 0$.

When $C=R$ one writes $\gdim_R(X)$ in place of $\gkdim{R}_R(X)$.  
When $X$ is a nonzero finite $R$-module, this is the G-dimension 
of~\cite{auslander:adgeteac,auslander:smt};  in general, it is the 
G-dimension of~\cite{christensen:gd,yassemi:gd}.

For each semidualizing complex, the $\mathrm{G}_C$-dimension is a 
refinement of the projective dimension and satisfies
an analogue of Theorem A.  Furthermore,   
a finite $R$-module of finite $\mathrm{G}_C$-dimension satisfies
an ``AB formula'';
see~\cite[(3.14),(3.15),(8.4)]{christensen:scatac}.

\medskip

\noindent \textbf{Theorem B.}  
\textit{Let $C$ be a semidualizing $R$-complex and $N$ a finite 
$R$-module.
\begin{enumerate}[\quad\rm(1)]
\item There is an inequality
$\gkdim{C}_R(N)\leq \pd_R(N)$
with equality when $\pd_R(N)$ is finite.
\item If $\gkdim{C}_R(N)$ is finite, then 
$\gkdim{C}_R(N)=\depth(R)-\depth_R(N)$.  
\item The following conditions are equivalent:
\begin{enumerate}[\quad\rm(i)]
\item  $C$ is dualizing;
\item  $\gkdim{C}_Rk<\infty$;
\item  $\gkdim{C}_RM<\infty$ for all finite $R$-modules $M$.
\end{enumerate}
\end{enumerate}
}

\medskip

Some important ring-theoretic properties are implied by 
the existence and behavior of dualizing modules.

\medskip

\noindent \textbf{Theorem C.}  
\textit{If $R$ possesses a dualizing module, then $R$ is 
Cohen-Macaulay.  If the $R$-module $R$ is dualizing, then $R$ is 
Gorenstein.}

\medskip

We shall make use of the following properties of 
$\mathrm{G}_C$-dimension contained 
in~\cite[(5.10),(6.5)]{christensen:scatac}.

\medskip

\noindent \textbf{Proposition D.}  
\textit{Let $C$ be a semidualizing complex, $M$ a finite $R$-module 
and $x\in\fm$ an $R$-regular element.  Set $R'=R/xR$.
\begin{enumerate}[\quad\rm(1)]
\item  The $R'$-complex $C'=C\lotimes_R R'$ is semidualizing.
\item \label{item2}
If $x$ is $M$-regular and $\gkdim{C}_R(M)<\infty$, then 
$\gkdim{C'}_{R'}(M/xM)=\gkdim{C}_R(M)$.
\end{enumerate}
}

\medskip

The following result is well-known for G-dimension.  We do not know of 
a reference for it in this generality, so we include a proof here.

\begin{lem} \label{lem1}
Let $C$ be a semidualizing $R$-complex and
$X\to Y\to Z\to\shift X$ a distinguished triangle in $\DD(R)$.  If
two of the complexes $X,Y,Z$ are $C$-reflexive, then so is the third.  
\end{lem}

It follows that, given an exact sequence of finite $R$-modules 
\[ 0\to M_r\to M_{r-1}\to\cdots\to M_1\to 0\]
with $\gkdim{C}_R(M_i)<\infty$ for $i\neq j$, 
one also has $\gkdim{C}_R(M_j)<\infty$.  Furthermore, it is 
straightforward to show that $\gkdim{C}_R(M_1\oplus M_2)<\infty$ 
implies $\gkdim{C}_R(M_i)<\infty$ for $i=1,2$.

\begin{proof}
Since $R$ is Noetherian, 
an analysis of the long exact sequence coming from the distinguished 
triangle $X\to Y\to Z\to\shift X$ shows that, if two of the complexes 
$X,Y,Z$ are homologically finite, so is the third.   Similarly for 
the distinguished triangle 
$\rhom_R(Z,C)\to\rhom_R(Y,C)\to\rhom_R(X,C)\to\shift\rhom_R(Z,C)$.  
The naturality of the biduality morphism gives rise to a commuting 
diagram
% \begin{diagram}[small]
% X & \rTo & \rhom_R(\rhom_R(X,C),C) \\
% \dTo & & \dTo \\
% Y & \rTo & \rhom_R(\rhom_R(Y,C),C) \\
% \dTo & & \dTo \\
% Z & \rTo & \rhom_R(\rhom_R(Z,C),C) \\
% \dTo & & \dTo \\
% \Sigma X & \rTo & \Sigma\rhom_R(\rhom_R(X,C),C) \\
% \end{diagram}
\[ \xymatrix{
X \ar[r] \ar[d] & \rhom_R(\rhom_R(X,C),C) \ar[d] \\
Y \ar[r] \ar[d] & \rhom_R(\rhom_R(Y,C),C) \ar[d] \\
Z \ar[r] \ar[d]  & \rhom_R(\rhom_R(Z,C),C) \ar[d] \\ 
\shift X \ar[r]  & \shift\rhom_R(\rhom_R(X,C),C)
} \]
where each column is a distinguished triangle.
Thus, if the biduality morphism is an isomorphism for 
two of the complexes $X,Y,Z$, then so is the third.
\end{proof}

Three homological dimensions of note have been introduced that 
characterize the complete intersection property like projective 
dimension does for regularity:  virtual projective dimension, 
complete intersection dimension, and lower complete intersection 
dimension.  Here we consider the last of these, as it is 
simultaneously the least restrictive and most flexible of the three.  
The interested reader is encouraged to consult the original 
sources~\cite{avramov:mofvpd,avramov:cid} 
for information on the other two.
The original treatment of lower complete intersection dimension is 
in~\cite{gerko:ohd}.  The (equivalent and slightly simpler) 
description we give here is from 
Sather-Wagstaff~\cite{sather:cidfc}.   

Quite simply, a finite $R$-module $M$ has finite lower complete 
intersection dimension, denoted $\lcidim_R(M)<\infty$, if it has finite 
G-dimension and the Betti numbers $\beta^n_R(M)$ are bounded above by 
a polynomial in $n$.  When $\lcidim_R(M)<\infty$ we set 
$\lcidim_R(M)=\gdim_R(M)$.  Basic properties of 
$\mathrm{CI}_{\ast}$-dimension are 
taken from~\cite{gerko:ohd} and summarized in the following.

\medskip

\noindent \textbf{Theorem E.}  
\textit{Let $M_1,\ldots,M_r$ be finite $R$-modules.
\begin{enumerate}[\quad\rm(1)]
\item  The following conditions on $R$ are equivalent.
\begin{enumerate}[\quad\rm(i)]
\item  $R$ is complete intersection;
\item  $\lcidim_Rk<\infty$;
\item  $\lcidim_RM<\infty$ for all finite $R$-modules $M$.
\end{enumerate}
\item Given an exact sequence 
\[ 0\to M_r\to M_{r-1}\to\cdots\to M_1\to 0\]
with $\lcidim_R(M_i)<\infty$ for $i\neq j$, 
one also has $\lcidim_R(M_j)<\infty$.  
\item If $\lcidim_R(M_1\oplus M_2)<\infty$, 
then $\lcidim_R(M_i)<\infty$ for $i=1,2$.
\item When $x\in\fm$ is $R$- and $M_1$-regular, set $R'=R/xR$
and $M_1'=M_1/xM_1$. If $\lcidim_R(M_1)<\infty$, then 
$\lcidim_{R'}(M_1')=\lcidim_R(M_1)$.
\end{enumerate}
}

\medskip

The final homological dimension we shall employ is the Gorenstein 
injective dimension of Enochs and Jenda~\cite{enochs1,enochs2}.  

An $R$-module $L$ is 
\emph{Gorenstein injective} if there exists an exact complex of injective 
$R$-modules $I$ such that the complex $\Hom_R(I,J)$ is 
exact for each injective $R$-module $J$ and 
$L\cong \coker(\partial_1^I)$.  
A \emph{Gorenstein injective resolution} of 
a module $M$ is a complex of Gorenstein injective modules $L$
with $L_i=0$ for all $i>0$, $\HH_i(L)=0$ for all $i\neq 0$, and 
$\HH_0(L)\cong M$.  The \emph{Gorenstein injective dimension} of $M$ 
is
\[ \gid_R(M)=\inf\{\sup\{n\in\mathbb{Z}\mid L_n\neq 0\} \mid
\text{$L$ a Gorenstein injective resolution of $M$}\}.\]

The finiteness of
Gorenstein injective dimension characterizes Gorenstein rings 
like the finiteness of injective dimension does for regular rings;  
see Khatami and Yassemi~\cite[(2.7)]{khatami}.

\medskip

\noindent \textbf{Theorem F.}  
\textit{The following conditions are equivalent:
\begin{enumerate}[\quad\rm(i)]
\item  $R$ is Gorenstein;
\item  $\gid_Rk<\infty$;
\item  $\gid_RM<\infty$ for all finite $R$-modules $M$.
\end{enumerate}
}

\medskip

The behavior of Gorenstein injective dimension with respect to exact 
sequences is described in the work of Holm~\cite{holm}.

\medskip

\noindent \textbf{Proposition G.}  
\textit{Let $M_1,\ldots,M_r$ be finite $R$-modules.
\begin{enumerate}[\quad\rm(1)]
\item Given an exact sequence 
\[ 0\to M_r\to M_{r-1}\to\cdots\to M_1\to 0\]
with $\gid_R(M_i)<\infty$ for $i\neq j$, 
one also has $\gid_R(M_j)<\infty$.  
\item If $\gid_R(M_1\oplus M_2)<\infty$, 
then $\gid_R(M_i)<\infty$ for $i=1,2$.
\end{enumerate}
}

\medskip

Gorenstein injective dimension also behaves well with respect to 
killing a regular element.  When $R$ possesses a dualizing complex, 
this is contained in the work of Christensen, Frankild, and 
Holm~\cite[(5.5)]{christensen:new}.

\begin{lem} \label{lem2}
Let $M$ be a finite $R$-module and 
$x\in\fm$ an $R$- and $M$-regular element.  Set $R'=R/xR$
and $M'=M/xM$. If $\gid_R(M)<\infty$, then 
$\gid_{R'}(M')<\infty$.
\end{lem}

\begin{proof}
For a Gorenstein injective module $L$, the $R'$-module $\Hom_R(R',L)$ 
is Gorenstein injective by~\cite[(3.1)]{enochs3}.
We argue as in~\cite[(3.10)]{christensen:new}.  Let $L$ be a bounded 
Gorenstein injective resolution of $M$, and $I$ an injective 
resolution of $M$ that is bounded above.  
Using 
Avramov and Foxby~\cite[(1.1.I),(1.4.I)]{avramov:hdouc} there is a 
quasi-isomorphism $L\xra{\simeq} I$.  
From~\cite[(3.10)]{christensen:new}, the induced homomorphism
$\Hom_R(R',L)\to\Hom_R(R',I)$ is a quasi-isomorphism as well.  
This morphism preserves the $R'$-structures, 
and $\Hom_R(R',L)$ is a bounded complex of Gorenstein injective 
modules. Thus, 
the isomorphism $\rhom_R(R',M)\simeq\rhom_R(R',L)$ in $\DD(R)$ is 
also an isomorphism in $\DD(R')$ and 
$\gid_{R'}(\rhom_R(R',M))<\infty$.  Since $x$ is $M$-regular, there 
are $R'$-isomorphisms
\[ M'\simeq R'\lotimes_R M\simeq \shift\rhom_R(R',M) \]
and it follows that $\gid_{R'}(M')<\infty$.
\end{proof}

\section{Characterizations of rings} \label{sec2}

The main result of this paper is the following.

\begin{thm} \label{thm1}
Let $C$ be a semidualizing $R$-complex and 
$M$ a finite $R$-module.
Set $d=\depth(R)$ and fix an integer $t$ 
between 0 and $d$.  If $\gkdim{C}_R(M/\fa M)$ is
finite for all ideals $\fa$ generated by an $R$-regular
sequence of length at most $d-t$ then either
$\gkdim{C}_R(M)\leq t$ or $C$ is a
dualizing complex.  
\end{thm}

\begin{proof}
Assume without loss of generality that $M$ is nonzero
and argue by induction on $d$.  If $d=0$, then $t=0$ and
$\gkdim{C}_R M=0$ by the AB-formula and we are done.  
Similar reasoning allows us 
to assume that $t<d$ for the rest of the proof.

Next, consider the case $d\geq 1$ and $\depth_R M=0$.  
Fix an $R$-regular element $x$
such that $x\in\fm\setminus\bigcup_{\fp\in S}\fp$,
where $S=\Ass M\backslash\{\fm\}$.  Since $M$ is Noetherian
there exists an integer $\alpha$ such that
$0:_Mx^{\alpha}=0:_Mx^{\alpha+i}$ for all $i\ge 0$.  The assumption
$\depth_R M=0$ yields $0:_M\fm\neq 0$ and hence
$0:_Mx^{\alpha}\neq 0$. Set $N=0:_Mx^{\alpha}$ and use the exact
sequence
\[ 0\to N\to M\stackrel{x^{\alpha}}{\rightarrow}M\to
M/x^{\alpha}M\to 0 \]
to deduce that $\gkdim{C}_RN<\infty$. 

The definition of $N$ implies that
every associated prime of $N$ contains $x$.  
On the other hand, $x$ is not in any prime in $S$, so the containment
$\Ass N\subseteq\Ass M$ gives
$\Ass N=\{\fm\}$.  Thus, there exists an integer $n\geq 1$ 
such that $\fm^nN=0$ and $\fm^{n-1}N\neq0$. 
If $n=1$ then $k$ is a direct summand
of $N$ so that $\gkdim{C}_Rk<\infty$ and $C$ is a dualizing
complex. If $n>1$ we can assume that $\fm^nN=0$ and $\fm^{n-1}N\neq
0$. It is straightforward to verify
that $N=\H^0_\fm(M)$, and therefore $\depth_R(M/N)>0$.
Fix $y\in\fm^{n-1}\setminus((\bigcup_{\fp\in S}\fp)
\cup(\bigcup_{\fp\in\Ass(M/N)}\fp)\cup\Ann N)$.  It follows that
$yN\neq 0$ and $y\fm N=0$.  Let $\beta$ be an integer such that
$0:_My^{\beta}=0:_My^{\beta+i}$ for all $i\ge 0$.  Arguing as above 
yields $0:_My^{\beta}=\H^0_\fm(M)=N$.

Applying the functor $R/y\otimes_R (-)$ to the exact sequence
\begin{equation}
0\to N\xra{\vf} M\to M/N\to 0 \tag{$\ast$} \label{eq1}
\end{equation}
induces the sequence
\begin{equation}
0\to R/y\otimes_R N\xra{\ol{\vf}} R/y\otimes_R M\to R/y\otimes_R (M/N)\to
0 \tag{$\ast\ast$} \label{eq2}
\end{equation}
which is exact because $\ker(\ol{\vf})=\Tor^R_1(R/y,M/N)=0$;  c.f., 
\cite[(16.5.i)]{matsumura:crt}.
From~\eqref{eq1} we have $\gkdim{C}_R(M/N)<\infty$ and hence
the exact sequence 
\[ 0\to M/N\xra{\cdot y} M/N\to R/y\otimes_R(M/N) \to 0 \]
yields
$\gkdim{C}_R(R/y\otimes_R(M/N))<\infty$.
Furthermore, 
$\gkdim{C}_R(R/y\otimes_RM)$ is finite and so~\eqref{eq2} gives
$\gkdim{C}_R(R/y\otimes_RN)<\infty$.  Similarly, the exact sequence
\[ 0\to yN\to N\to N/yN\to 0 \]
yields $\gkdim{C}_R(yN)<\infty$. Since
$\fm(yN)=0$, the $R$-module $k$ is a direct summand of $yN$
and so $\gkdim{C}_R(k)<\infty$. Therefore $C$ is dualizing.

Now suppose that $\depth R\ge 1$ and 
$\depth_RM\geq 1$. 
Fix an element
$z\in \fm$ that is both $R$- and $M$-regular,
and set $\bar{R}=R/xR$ and 
$\bar{M}=M/xM$. The $\bar{R}$-complex $\bar{C}=C\lotimes_R \bar{R}$ 
is semidualizing, and it follows easily from Proposition D that 
$\gkdim{\bar{C}}_{\bar{R}}(\bar{M}/\mathfrak{b}\bar{M})$ is finite for 
each ideal $\mathfrak{b}$ generated by an $\bar{R}$-regular sequence 
of length at most $d-1-t$.  Thus, by induction, either 
$\bar{C}$ is dualizing for $\bar{R}$ or
$\gkdim{\bar{C}}_{\bar{R}}(\bar{M})\leq t$.  If 
$\bar{C}$ is dualizing for $\bar{R}$, then
$C$ is dualizing for $R$.  If 
$\gkdim{\bar{C}}_{\bar{R}}(\bar{M})\leq t$, then
$\gkdim{C}_{R}(M)\leq t$.
\end{proof}

Applying Theorem~\ref{thm1} to a semidualizing \emph{module}
yields a
criterion for the Cohen-Macaulay property
that is parallel to Theorem A2.

\begin{cor} \label{cor1}
Let $C$ be a semidualizing $R$-module and 
$M$ a finite $R$-module.
Set $d=\depth(R)$ and fix an integer $t$ 
between 0 and $d$.  If $\gkdim{C}_R(M/\fa M)$ is
finite for all ideals $\fa$ generated by an $R$-regular
sequence of length at most $d-t$ then either
$\gkdim{C}_R(M)\leq t$ or $R$ is Cohen-Macaulay.
\end{cor}

Using $C=R$ gives a criterion for the Gorenstein property.
A similar result in terms of upper G-dimension follows immediately 
from this one.

\begin{cor} \label{cor2}
Let $M$ be a finite $R$-module.
Set $d=\depth(R)$ and fix an integer $t$ 
between 0 and $d$.  If $\gdim_R(M/\fa M)$ is
finite for all ideals $\fa$ generated by an $R$-regular
sequence of length at most $d-t$ then either
$\gdim_R(M)\leq t$ or $R$ is Gorenstein.
\end{cor}

Modifying the proof of Theorem~\ref{thm1} appropriately, one obtains 
analogous results using virtual projective dimension, complete 
intersection dimension, or lower complete intersection dimension to 
describe complete intersection rings.  We state here the 
version for lower complete intersection dimension, as the others can 
be derived immediately from it.  

\begin{cor} \label{cor4}
Let $M$ be a finite $R$-module.
Set $d=\depth(R)$ and fix an integer $t$ 
between 0 and $d$.  If $\lcidim_R(M/\fa M)$ is
finite for all ideals $\fa$ generated by an $R$-regular
sequence of length at most $d-t$ then either
$\lcidim_R(M)\leq t$ or $R$ is complete intersection.
\end{cor}

A similar modification of the proof of Theorem~\ref{thm1}
produces a criterion for the regularity property that generalizes 
Theorem A2.  

\begin{cor} \label{cor3}
Let $M$ be a finite $R$-module.
Set $d=\depth(R)$ and fix an integer $t$ 
between 0 and $d$.  If $\pd_R(M/\fa M)$ is
finite for all ideals $\fa$ generated by an $R$-regular
sequence of length at most $d-t$ then either
$\pd_R(M)\leq t$ or $R$ is regular.
\end{cor}

The analogous statements for injective dimension 
and Gorenstein injective dimension 
differ slightly from the previous results because 
of the 
Bass formula.  The proof is the same, though, 
modulo easy arguments to deal with the first cases. 

\begin{cor} \label{cor6}
Let $M$ be a finite $R$-module.
Set $d=\depth(R)$ and fix an integer $t$ 
between 0 and $d$.  If $\id_R(M/\fa M)$ is
finite for all ideals $\fa$ generated by an $R$-regular
sequence of length at most $d-t$ then either
$\depth_R M\geq d-t$ or $R$ is regular.
\end{cor}

\begin{cor} \label{cor5}
Let $M$ be a finite $R$-module.
Set $d=\depth(R)$ and fix an integer $t$ 
between 0 and $d$.  If $\gid_R(M/\fa M)$ is
finite for all ideals $\fa$ generated by an $R$-regular
sequence of length at most $d-t$ then either
$\depth_R M\geq d-t$ or $R$ is Gorenstein.
\end{cor}

Our final variation on this theme has a similar proof, but considers
ideals generated by parts of system of parameters instead of regular
sequences.

\begin{cor} \label{cor06}
Let $C$ be a semidualizing $R$-module.  The
following conditions on $R$ are equivalent.
\begin{enumerate}[\quad\rm(i)]
\item $R$ is Cohen-Macaulay.
\item There exists a finite $R$-module $M$ such that for 
every ideal $\fa$ generated by part of a
system
of parameters for $R$, one has $\gkdim{C}dim_R(M/\fa M)$
finite.
\item For every ideal $\fa$ generated by part of a
system
of parameters for $R$, one has $\gkdim{C}dim_R(M/\fa M)$
finite.
\end{enumerate}
\end{cor}

It would be interesting to know whether the parallel result using 
Cohen-Macaulay dimension to characterize Cohen-Macaulay rings holds.  
The only obstruction is the 
current lack of understanding of the behavior 
of this homological dimension with respect to exact sequences.

\section*{Acknowledgments} 

The first and third authors were visiting the Abdus Salam 
International
Centre for Theoretical Physics (ICTP)
during the preparation of this paper. They would like to thank
the ICTP for its hospitality during their stay there.

\providecommand{\bysame}{\leavevmode\hbox to3em{\hrulefill}\thinspace}

\end{document}